\documentclass[12pt]{article}

\usepackage[utf8]{inputenc}
\usepackage[english]{babel}
 \usepackage[T1,T2A]{fontenc}
\usepackage{amssymb,amsmath,amsthm}
\usepackage{graphicx}
\usepackage{color}
\usepackage{ulem}
\usepackage{hyperref}
\usepackage{extarrows}
\usepackage{floatflt}
\usepackage{caption}
\usepackage{framed}
\usepackage{amscd,amsfonts,amssymb,amsmath,latexsym,array,hhline,xcolor,graphicx}
\usepackage[justification=centering]{caption}

\titlepage
\newtheorem*{thm*}{Теорема}

\newtheorem*{lem*}{Лемма}

\newtheorem*{prop*}{Утверждение}

\newtheorem*{crl*}{Следствие}
\newtheorem*{def*}{Определение}

\newcommand\bean{\begin{eqnarray*}}
\newcommand\eean{\end{eqnarray*}}
\newcommand\beq{\begin{equation}}
\newcommand\eeq{\end{equation}}

\begin{document}
 \begin{Huge}
\centerline{On discrete part of Dirichlet spectrum}
\end{Huge}
\begin{Large}
\centerline{by Sergei Pitcyn}
 \end{Large}

\section{Introduction}

In this paper we consider rational approximations to an irrational number ${\alpha \in \mathbb{R}}$.
Quite recently 
Han\v{c}l \cite{Hc0,Hc} obtained a result which improves on approximations to real numbers which correspond to the discrete part of Lagrange spectrum.
In the present paper we prove a similar result related to the discrete part of Dirichlet spectrum.
The paper is organized as follows.  In Subsection \ref{1.1} of our Introduction we recall classical results related to Lagrange spectrum and formulate the improvement by 
Han\v{c}l, and 
in Subsection \ref{1.2} we describe well-known results concerning Dirichlet spectrum. Our main new results are presented in Section \ref{2.}. We formulate two theorems. Theorem 1 deals just with the minimal element of Dirichlet spectrum, while Theorem 2 takes into consideration the whole discrete part of Dirichlet spectrum, which was discovered by Morimoto in 1930th. In Section \ref{3.} we give a proof of Theorem 1 and in Sections \ref{4.}, \ref{5.} we give a proof of Theorem 2.

\subsection{Lagrange spectrum}\label{1.1}

  Let us consider the irrationality measure function
$$
\psi_\alpha (t) = \min_{q\in \mathbb{Z}_+: 1\le q \le t }
||q\alpha||,\,\,\,\, \text{where}\,\,\,\, ||\xi || = \min_{a\in \mathbb{Z}}
|\xi - a|.
$$
Let 
$$ \alpha = [a_0;a_1,\,\cdots,a_n,\,\cdots], \,\,\,\, a_0 \in \mathbb{Z}, \,\,\,\, a_j \in \mathbb{Z}_+, j=1,2,3,\cdots
$$
be the representation of $\alpha$ as a continued fraction and $ \frac{p_n}{q_n} =[a_0;a_1,\,\cdots,a_n] $ be its $n$-th convergent. 
Numbers $ \alpha = [a_0;a_1,\,\cdots,a_n,\,\cdots]$ and 
$\beta  = [b_0;b_1,\,\cdots,b_m,\,\cdots]
 $
 are called {\it equivalent} if there exist $n_0$ and $m_0$
such that $a_{n_0+s}= b_{m_0+s}$ for all $ s =1,2,3,\cdots$.
If $\alpha$ and $\beta$ are equivalent, we write $ \alpha \sim \beta$.

By the Lagrange theorem we have the equality
$$
\psi_\alpha (t) = || q_n\alpha||\,\,\,\,
\text{for}\,\,\,\,\,
q_t \le t <q_{n+1}.
$$

We consider the Lagrange constant
$$
\lambda(\alpha) = \left(\liminf_{t\to\infty} t\cdot \psi_\alpha (t)\right)^{-1} =
 \left(\liminf_{n\to\infty} q_n\cdot ||q_n\alpha||\right)^{-1}\!.
$$

For equivalent numbers $\alpha \sim \beta$ one has 
$\lambda(\alpha) = \lambda(\beta)$.
The set $\mathbb{L}$ of all possible values of $\lambda(\alpha)$ is known as the
{\it Lagrange spectrum}. A famous theorem by Markoff (see 
books \cite{Cas,Cus} for detailed exposition and discussions) claims that the set $\mathbb{L}\cap \left[0, 3\right)$
is countable and consists of increasing numbers 
$$
L_1   = {\sqrt{5}}<
L_2  = {\sqrt{8}}<
L_3  = \ {\sqrt{221}/5}<
\cdots
<
L_j = L(M_j) =  {\sqrt{9- \frac{4}{M_j^2}}}
,
\,\cdots
,
$$
where $M_j, j=1,2,3,\cdots$ run over all the so-called {\it Markoff numbers}, and there exist certain quadratic irrationalities $\gamma_j$ such that the equality $\lambda(\alpha) =L_j$ holds if and only if $\alpha \sim \gamma_j$.

Further consideration 
(see Theorem III from Ch. II \cite{Cas})
shows  that the following statement is valid.

\vskip+0.3cm
{\bf Proposition A.}
{\it Inequality 
\begin{equation}\label{wq1}
||q\alpha|| < \frac{1}{L_m q}
\end{equation}
has infinitely many solutions in positive integers $q$
(or equivalently, 
the set of
numbers  $t$ satisfying
$
\psi_\alpha (t) <\frac{1}{L_m t}
$
is unbounded)
if and only if $\alpha\not\sim \gamma_j$ for all $j\le m-1$.}

\vskip+0.3cm

We should note that the famous theorem by Hurwitz, which claims that for any irrational $\alpha$
there exist infinitely many $q$ with ${||q\alpha||<\frac{1}{\sqrt{5}q}}$, may be considered as a particular case with $m=0$.

\vskip+0.3cm
In \cite{Hc}
Han\v{c}l gave an improved version of the previous statement.
One of his formulations is as follows.

\hskip+0.3cm
{\bf Theorem A.}
{\it 
  Consider function
$$
g_m(x) = \frac{2}{L_m\left(1+\sqrt{1+\frac{4}{L_m^2x^2}}\right) x }  .
$$
Assume  that  $\alpha\not\sim \gamma_j$ for all $j\le m-1$. Then inequality 
\begin{equation}\label{wq2}
||q\alpha|| \le  g_m(q)
\end{equation}
has infinitely many solutions in positive integers $q$.}
\hskip+0.3cm

Here we should note that the right-hand side of (\ref{wq2}) is smaller than that of (\ref{wq1}), that is $g_m(x)< \frac{1}{L_m x}$ for all $x \ge 1$.

\vskip+0.3cm
  Han\v{c}l proved that the result of Theorem A is optimal in a certain special sense. An earlier result by Han\v{c}l  \cite{Hc0} deals with just the minimal element $\sqrt{5}$ in $\mathbb{L}$.

\subsection{Dirichlet spectrum}\label{1.2}

In this subsection we describe the main results concerning the discrete part of the Dirichlet spectrum.
The Dirichlet constant for an irrational $\alpha$ is defined as 
$$
D(\alpha) = \limsup_{t\to\infty} t\cdot \psi_\alpha (t) =
\limsup_{n\to\infty} q_{n+1}\cdot ||q_n\alpha||.
$$

The set of all possible values of $D(\alpha)$ is known as the 
{\it Dirichlet spectrum} $\mathbb{D}$. It is clear that if $\alpha\sim\beta$, then $D(\alpha) = D(\beta)$.

Consider numbers
\begin{equation}\label{wq3}
\alpha_0 = [\overline{1}] = \frac{1+\sqrt{5}}{2},\,\,\,
\alpha_k =[\overline{\underbrace{1;\,\cdots,1}_{2k-1},2}] = [\overline{1_{2k-1},2}],\,\,\, k =1,2,3,\cdots
\end{equation}
(here and in the sequel we use notation 
$1_s$ for $s$ consequtive partial quotients
$ \underbrace{1,\cdots, 1}_s$),
defined by their continued fractions, and the values 
$$ D_0 = D(\alpha_0)
=\frac{1}{2}+\frac{1}{2\sqrt{5}}< D_1 = D(\alpha_1)<\cdots< D_k = D(\alpha_k)<\cdots.
$$

These numbers were discovered by Morimoto \cite{Morimo}. In particular
\begin{equation}\label{Dkk}
D_k = \frac{2\alpha_k+1}{2\alpha_k+2},\,\,\,\, k \ge 1.
\end{equation}

Lesca in his thesis \cite{Lesca} gave a complete proof of the following equality:
$$\mathbb{D}\cap \left[\frac{1}{2}+\frac{1}{2\sqrt{5}},\frac{1+\sqrt{5}}{4}\right)
= \{D_j:\,j =0,1,2,3,\cdots \}.
$$ 

Moreover, he showed that equality $D(\alpha) =D_j$ holds if and only if $ \alpha \sim \alpha_j$.

\vskip+0.3cm
In particular, as a corollary of this result we can get the following statement.

\newpage
{\bf Proposition B.}
 
\noindent
{\bf a}) {\it If
$\alpha$ satisfies $\alpha \not\sim \alpha_j$ for all $j$ from the interval 
$0\le j \le m-1$, then  
$$
\limsup_{t\to \infty} t\cdot \psi_\alpha (t) \ge D_m
$$
holds;}
\vskip+0.3cm
\noindent
{\bf b}) {\it there exist $\alpha$, 
satisfying $\alpha \not\sim \alpha_j$
for all $j$ from the interval 
$0\le j \le m-1$,
such that 
$$
\limsup_{t\to \infty} t\cdot \psi_\alpha (t) = D_m;
$$}
\vskip+0.3cm
\noindent
{\bf c}) {\it all $\alpha$'s satisfying {\rm {\bf b})} are equivalent to $\alpha_m$.}
 \vskip+0.3cm

 The case $m=0$ recovers a famous result by Szekeres \cite{sec10}: for any irrational $\alpha$ there exist
 infinitely many $n$ such that
 \begin{equation}\label{sec1}
 \lim_{t\to q_{n+1}-} t\cdot \psi_\alpha (t) =
 q_{n+1}||q_n\alpha||> D_0 = 
 \frac{1}{2}+\frac{1}{2\sqrt{5}}.
 \end{equation}

A brief description of some other results related to Dirichlet spectrum one can find, for example, in relatively recent paper \cite{A}.

\section{Main results}\label{2.}

In Subsection \ref{1.1} we recalled the results related to the discrete part of the Lagrange spectrum and formulated relatively recent improvements from \cite{Hc}. 
In the present paper we,
following   Han\v{c}l's ideas in the case of the Lagrange spectrum,
deduce an improvement of Morimoto-Lesca's results about the discrete part of the Dirichlet spectrum, discussed in Subsection \ref{1.2}.
In Theorem 1 below we improve on Szekeres' result \ref{sec1}. Theorem 2 deals with an improvement of general Proposition B for all elements of the discrete part of $\mathbb{D}$.
\vskip+0.3cm

Let us define the function
$$
f_0(x) = \frac{\sqrt{5} + 1}{\sqrt{5} + \sqrt{5 - \frac{4}{x^{2}}}}.
$$
We should note that 
 \begin{equation}\label{sec10}
f_0(x) > D_0=\frac{1}{2}+\frac{1}{2\sqrt{5}}  
 \end{equation}
 for all positive $x$.

\textbf{Theorem 1.} 

1. \textit{Let $\alpha$ be an irrational number. Then there exists an increasing sequence of positive integers $q_n$ such that}
 \begin{equation}\label{sec11}
 \lim_{t\to q_{n}-} t\cdot \psi_\alpha (t) \ge
 f_0 (q_n).
 \end{equation}

2. \textit{Two conditions}

a) $ t\psi_{\alpha}(t)<f_{0}(t)$ {\it for all $t$ large enough};

b) {\it there exists a sequence of positive integers $q_n$ such that}
  $\lim\limits_{t \to q_{n} - } \frac{t\psi_{\alpha}(t)}{f_{0}(t)}=1$;

{\it are satisfied simultaneously if and only if 
  $\alpha = \pm\alpha_0 +C$ with $ C\in \mathbb{Z}$}.

\vskip+0.3cm

Remark 1. Inequality (\ref{sec11}) is stronger than (\ref{sec1}) because of (\ref{sec10}).

\vskip+0.3cm

Remark 2. Certainly, the sequence $q_n$ from Theorem 1 forms a subsequence of the sequence of denominators of convergents of $\alpha$.

 \vskip+0.3cm

In addition, in this work a similar strengthened result is proven for all the numbers at which the discrete part of the Dirichlet spectrum is reached (that is, for the numbers $\alpha \sim \alpha_{k} = [\overline{1_{2k-1}, 2}]$ for some  $k \in \mathbb N$).

Let us 
define numbers
\begin{equation}\label{betamain}
\beta_{k} = [0;1_{k-1}, \overline{2, 1_{2k-1}}], \,\,\,\beta^{(1)}_{k} = [0;1_{k}, \overline{2, 1_{2k-1}}],\,\,\,\beta^{(2)}_{k} = [0;1_{k-2}, \overline{2, 1_{2k-1}}]
\end{equation}
and
introduce the function
$$
f_k(x) = 
\begin{cases}
  \frac{2D_{k}}{1 + \sqrt{ 1 - \frac{2 \alpha_{k}}{(2 \beta_{k} + 1)( \alpha_{k} + 1) x^{2}}}} , k \equiv 1  \pmod{2}, \cr
  \frac{2D_{k}}{1 + \sqrt{ 1 - \frac{2 \alpha_{k}}{(\beta^{(1)}_{k} + \beta^{(2)}_{k} + 1)( \alpha_{k} + 1) x^{2}}}}, k \equiv 0 \pmod{2},
\end{cases}
$$
where $\alpha_k$ and $D_k$ are defined in (\ref{wq3}) and (\ref{Dkk}) respectively.
It is clear that
$$
f_k(x)> D_k\,\,\,\,\,\, \forall\, x >1,\,\,\, \forall\, k.
$$
Then the following statement holds.

\textbf{Theorem 2.}

1. \textit{Let $\alpha \not\sim \alpha_j$ for all $j$ from the interval $0\le j \le k-1$. Then there exists an increasing sequence of positive numbers $q_n$ such that
\begin{equation}\label{sec12}
 \lim_{t\to q_{n}-} t\cdot \psi_\alpha (t) \ge
 f_k (q_n).
 \end{equation}}

2. \textit{ 
Let
$$
\alpha= 
\begin{cases}
    \beta_k\,\,\, \text{if} \,\, k \,\,\text{is odd},\cr
    \beta_k^{(1)}\,\,\text{or}\,\, \beta_k^{(2)}\,\,\, \text{if} \,\, k \,\,\text{is even}.
\end{cases}
$$
Then}

a) $ t\psi_{\alpha}(t)<f_{k}(t)$ {\it for all $t$ large enough};

b) {\it there exists a sequence of positive integers $q_n$ such that} $\lim\limits_{t \to q_{n} - } \frac{t\psi_{\alpha}(t)}{f_{k}(t)}=1$.

\section{Proof of Theorem 1} \label{3.}

By Perron's formula, for $\alpha_0$ we have
\begin{equation}\label{pero}
q_{n+1}||q_{n} \alpha_0||  = \frac{1}{1 + [0;1_{n+1}] \cdot [0; \overline{1}]} = \frac{1}{1 + \frac{q_n}{q_{n+1}} \cdot \frac{1}{\alpha_0}}.
\end{equation}
It is clear that $\lim\limits_{n \to \infty} [0;1_{n+1}] = [0;\overline{1}] = \frac{1}{\alpha_0}$ and $\frac{q_{n}}{q_{n+1}}$ is just the $n+1$-th convergent of $\frac{1}{\alpha_0}.$ Then
$$\frac{1}{\alpha_0} - \frac{q_{n}}{q_{n+1}} = \frac{(-1)^{n+1}}{q_{n+1}^{2}([0;1_{n+1}] + [\overline{1}])} = \frac{(-1)^{n+1}}{q_{n+1}^{2}(\frac{q_n}{q_{n+1}} + \alpha_0)}$$
and so
$$\frac{q_{n}}{q_{n+1}} = \frac{-1 + \sqrt{5 + \frac{4  (-1)^{n+1}}{q_{n+1}^{2}}}}{2}$$
(in particular, the expression in square root is a full square of a rational number).
We substitute this into (\ref{pero}) to obtain
$$
\lim_{t\to q_{n+1}-} t\psi_{\alpha_0} (t)
=
q_{n+1}||q_{n} \alpha_0|| = \frac{\sqrt{5} + 1}{\sqrt{5} + \sqrt{5 + \frac{4(-1)^{n+1}}{q_{n+1}^{2}}}}
\le \frac{\sqrt{5} + 1}{\sqrt{5} + \sqrt{5 - \frac{4}{q_{n+1}^{2}}}}
$$
for all values of $n$. Hence
$$\psi_{\alpha_0}(t) <  \frac{\sqrt{5} + 1}{t\left(\sqrt{5} + \sqrt{5 - \frac{4}{t^{2}}}\right)} 
\,\,\,\,\,
\text{for all}\,\,\, t,
$$ 
and we have the statement 2 a) proven for $\alpha=\alpha_0$. Moreover, if we take $n=2\nu$ to be even, then 
$$
\lim\limits_{t \to q_{2\nu} - } \frac{t\psi_{\alpha_0}(t)}{f(t)}=1,$$
and we show that 2 b) is also valid for $\alpha = \alpha_0$.

\vskip+0.3cm
Now we prove statement 1.
For  $\alpha \nsim \alpha_0$
everything is clear because of Proposition A ($m=1$).
Therefore, we should consider only the case $\alpha \sim \alpha_0$.
Let  $\alpha = [a_{0};a_{1}, \cdots, a_{s}, \overline{1}].$ 
We consider  the denominator 
$q_{i}$ of the $i$-th convergent of $\alpha.$ Define $r_{s} = [a_{s};a_{s-1}, \cdots, a_{1}]$. Then
for $ n\ge s$ we have
\begin{equation}\label{pero1}
q_{n+1}||q_{n} \alpha|| = \frac{1}{1 + [0;1_{n+1-s}, r_{s}] \cdot [0; \overline{1}]} = \frac{1}{1 + \frac{q_n}{q_{n+1}} \cdot \frac{1}{\alpha_0}}.
\end{equation}
Certainly,
 $\lim\limits_{n \to \infty} [0;1_{n+1-s}, r_{s}] = [0;\overline{1}] = \frac{1}{\alpha_0}$,
 but now  
 $\frac{q_n}{q_{n+1}}$ 
 is not a convergent of $\frac{1}{\alpha_0}$, unless the period starts at the very beginning of the fraction.

  Now we deal with the difference
 $$\frac{1}{\alpha_0} - \frac{q_n}{q_{n+1}} = \frac{F_{n+1-s} \alpha_0 + F_{n-s}}{F_{n+2-s} \alpha_0 + F_{n+1-s}} - \frac{F_{n+1-s} r_{s} + F_{n-s}}{F_{n+2-s} r_{s} + F_{n+1-s}}
  = \frac{(-1)^{n-s}(\alpha_0 - r_{s})}{(F_{n+2-s} \alpha_0 + F_{n+1-s})(F_{n+2-s} r_{s} + F_{n+1-s})},$$

  where 
$F_i$ stands for $i$-th  Fibonacci number.
  
We express the value of the ratio $\frac{q_n}{q_{n+1}}$ from this equality and substitute it into (\ref{pero1}). Then we get 
$$q_{n+1}||q_{n} \alpha|| =  \frac{1}{1 + (\frac{1}{\alpha_0} +  \frac{(-1)^{n-s+1}(\alpha_0 - r_{s})}{(F_{n+2-s} \alpha_0 + F_{n+1-s})(F_{n+2-s} r_{s} + F_{n+1-s})})\cdot \frac{1}{\alpha_0}}=
\frac{\sqrt{5}+1}{\sqrt{5}+(\sqrt{5} + \phi_{n}(r_s))},
$$
where
$$
\phi_{n} (t) = \frac{(-1)^{n-s+1}2(\alpha_0 - t)}{(F_{n+2-s} \alpha_0 + F_{n+1-s})(F_{n+2-s} t + F_{n+1-s})}.
$$

 We finish the proof similarly to the argument from \cite{Hc0} (see cases 1 and 2 from the proof of Theorem 4.5 from \cite{Hc0}).

 In the case where $a_s \geqslant 3$ or $s>1$ we have the inequality
   $\frac{|\phi_{n}(t)|}{2}t^{2} > \frac{1}{\sqrt{5}} + \varepsilon$ 
   for any 
   $\varepsilon>0$ for $n$ large enough.
   This is just the inequality from case 1 from \cite{Hc0}.
   So if the difference  $n-s$ is even, we get    $\sqrt{5} + \phi_n(t) < \sqrt{5-\frac{4}{t^2}}$. This means that in our case 
 $$\lim_{t\to q_{n}-} t\cdot \psi_\alpha (t) >
 f_0 (q_n)$$
for all $\nu$ large enough.
 
The remaining case is $\alpha = [a_{0}; 2, \overline{1}]$.
This corresponds to case 2 from \cite{Hc0}. Following the argument from \cite{Hc0}, we see that 
  $$q_{n+1}||q_{n} \alpha|| = \frac{\sqrt{5}+1}{\sqrt{5} + \sqrt{5 + \frac{4 (-1)^{n}}{q_{n+1}^{2}}}}
  \ge f_0(q_{n+1})
  \,\,\,\, ( \, = f_0(q_{n+1})\,\, \text{for odd}\,\, n)
  .$$
  
Statement 1 in proven.
Moreover, we see that the number
$[a_0;2, \overline{1}]= a_0+2-\alpha_0$ also satisfies the conditions from Statement 2 of Theorem 1.
So  the consideration of the cases above shows that  the only numbers satisfying statement 2 are of the form 
  $\alpha_{0} + C$ and $C-\alpha_0$,  $C\in \mathbb{Z}$.

\section{To  Theorem 2: analysis of numbers $\alpha_k^{\pm 1}$} \label{4.}

The case $\alpha \sim \alpha_j $ for some $ 0\le j \le k-1$ is clear.
If
$\alpha \not\sim \alpha_j $ for all $ 0\le j \le k-1$ and in addition $\alpha \not\sim \alpha_k$ everything follows from Proposition B.
So we need only to consider the case when $\alpha\sim \alpha_k$.
We should work with numbers $\alpha_k, \frac{1}{\alpha_k}$
and all the numbers of the form $ [a_0;a_1,...,a_s,\alpha_k]$ with integers $a_j$.

\vskip+0.3cm
We start our proof by considering number  $\alpha_{k} = [\overline{ 1_{2k-1}, 2}]$ itself. From  \cite{Lesca} it is known that the values of $q_{n+1}||q_{n} \alpha_{k}||$ for sufficiently large $n$ {  admit local maximum } at ${n = n_{m} = 2km - 3}$ or ${n = n_{m} = 2km - 2}$ for sufficiently large natural $m$. To be more precise,
\begin{equation}\label{precise}
\limsup\limits_{n \to \infty} q_{n+1}||q_{n} \alpha_{k}|| = \lim\limits_{m \to \infty}  q_{2km-1}||q_{2km-2} \alpha_{k}|| =\lim\limits_{m \to \infty}  q_{2km-2}||q_{2km-3} \alpha_{k}||.
\end{equation}
At first, let ${n  = 2km - 3}$. Then
\begin{equation}\label{pero1}
q_{n+1}||q_{n} \alpha_{k}|| = \frac{1}{1 + [0;1_{2k-1},2, \,\cdots, 2, 1_{2k-2}] \cdot [0; \overline {2, 1_{2k-1}}]}.
\end{equation}

 Here $[0;1_{2k-1},2, \,\cdots, 2, 1_{2k-2}] = \frac{q_n}{q_{n+1}}$, where $q_{i}$ is the denominator of the $i$-th convergent of the number $\alpha$.

Obviously, $\lim\limits_{n \to \infty} \frac{q_n}{q_{n+1}} = [0; \overline{1_{2k-1}, 2}] = \frac{1}{\alpha_{k}}.$ Since in the present case the number $n+1$ is even and $\frac{q_{n}}{q_{n+1}}$ is a convergent of $\frac{1}{\alpha_{k}}$, { we have  equality}
\begin{equation}\label{re10}
\frac{1}{\alpha_{k}} - \frac{q_{n}}{q_{n+1}} = \frac{1}{q_{n+1}^{2}([0;1_{2k-2},2, \,\cdots, 2, 1_{2k-1}] + [1; \overline{2, 1_{2k-1}}])} = \frac{1}{q_{n+1}^{2}(\frac{p_{n+1}^{*}}{q_{n+1}} + \beta)}.
\end{equation}

We defined here $p_{n+1}^*$ and $\beta$  to satisfy the equalities $$[0;1_{2k-2}, \underbrace{2, 1_{2k-1},\,\,\cdots\,\, ,2, 1_{2k-1}}_{m-1\, \text{times}}] = \frac{p_{n+1}^{*}}{q_{n+1}}\,\,\,\text{and}\,\,\, [1; \overline{2, 1_{2k-1}}] = \beta = \frac{3 \alpha_{k} + 1}{2 \alpha_{k} + 1}.$$ 
 So ${\lim\limits_{n \to \infty} \frac{p_{n+1}^{*}}{q_{n+1}} = [0;1_{2k-2}, \overline{2, 1_{2k-1}}] = \alpha_{k} - 1}$. Applying Perron's formula to approximation of the number $\alpha_k -1$ by convergent $\frac{p_{n+1}^*}{q_{n+1}}$ we obtain
\begin{equation}\label{re1}
\alpha_{k} - 1 - \frac{p_{n+1}^{*}}{q_{n+1}} = \frac{1}{q_{n+1}^{2}([0;1_{2k-1}, 2, \,\cdots, 2, 1_{2k-2}] + [2; \overline{1_{2k-1}, 2}])} = \frac{1}{q_{n+1}^{2}(\frac{q_n}{q_{n+1}} + \gamma)},
\end{equation}
where 
\begin{equation}\label{re2}
\gamma = [\overline{2; 1_{2k-1}}] = \frac{2 \alpha_{k} + 1}{\alpha_{k}}.
\end{equation}
Then we substitute the value of $\frac{p_{n+1}^{*}}{q_{n+1}}$ from (\ref{re1}) into (\ref{re10}) and obtain
$$\frac{1}{\alpha_{k}} - \frac{q_n}{q_{n+1}} = \frac{1}{q_{n+1}^{2}( \beta + \alpha_{k} - 1 - \frac{1}{q_{n+1}^{2}(\frac{q_{n}}{q_{n+1}} + \gamma)})},$$

hence
$$ \frac{q_n}{q_{n+1}} = -1 + \frac{\alpha_{k} + 1}{\alpha_{k}} \sqrt{1 - \frac{2 \alpha_{k} + 1}{(\alpha_{k} + 1)^{2} q_{n+1}^{2}}}.$$

Now we substitute the 
this expression for $\frac{q_n}{q_{n+1}}$
 into (\ref{pero1}) and get equality
$$q_{n+1}||q_{n} \alpha_{k}|| = \frac{2 \alpha_{k} + 1}{(\alpha_{k} + 1)(1 + \sqrt{1 - \frac{2 \alpha_{k} + 1}{(\alpha_{k} + 1)^{2} q_{n+1}^{2}}})} = \frac{2 D_{k}}{1 + \sqrt{1 - \frac{4 D_{k}(1 - D_{k})}{q_{n+1}^{2}}}} = g_{k}(q_{n+1}),$$
where  we  have defined
\begin{equation}\label{re3}
g_{k}(x) = \frac{2 D_{k}}{1 + \sqrt{1 - \frac{4 D_{k}(1 - D_{k})}{x^{2}}}} >D_k.
\end{equation}
So in the case $ n = 2km-3$ we see that $q_{n+1}||q_{n} \alpha_{k}||>D_k$.

Now let $n = 2km - 2$ and $m$ be a sufficiently large natural number. Then
$$q_{n+1}||q_{n} \alpha_{k}|| = \frac{1}{1 + [0;2, 1_{2k-1}, \,\cdots, 2, 1_{2k-2}] \cdot [0; \overline {1_{2k-1}, 2} ]} = \frac{1}{1 + \frac{q_n}{q_{n+1}} \cdot \frac{1}{\alpha_{k}}},$$
and $\lim\limits_{n \to \infty} \frac{q_n}{q_{n+1}} = [0; \overline {2, 1_{2k-1}}] = \frac{1}{\gamma}$,
where $\gamma$ is defined in (\ref{re2}).
But now ${n+1}$ is odd and so we have $\frac{1}{\gamma} - \frac{q_n}{q_{n+1}} < 0$. Then $\frac{q_n}{q_{n+1}} > \frac{1}{\gamma}$, therefore $q_{n+1}||q_{n} \alpha_{k}|| < D_{k}$. So this case does not lead to an extremal value of $t\psi_{\alpha_k}(t)$.

To summarize, we proved the following

\textbf{Proposition 1.} \textit{
Suppose that $g_k(x)$ defined in (\ref{re3}).
Let $k>0$. Then the following two conditions are satisfied:}

a) $ t\psi_{\alpha_{k}}(t)<g_{k}(t)$ {\it for  
all $t$ large enough};

b) {\it there exists a sequence of positive integers $q_n$ such that}
  $\lim\limits_{t \to q_{n} - } \frac{t\psi_{\alpha_{k}}(t)}{g_{k}(t)}=1$.
\vskip+0.3cm

\textbf{Remark 1.} Proposition 1 remains valid if instead of number $\alpha_k$ we consider the number
 $\frac{1}{\alpha_{k}} = [0; \overline{ 1_{2k-1}, 2}]$,
 that is conditions  a) and b) are satisfied for irrationality measure function
  $\psi_{\frac{1}{\alpha_k}}(x)$. The calculations are quite similar. The only difference is that we need to change the parity and the extremal value which corresponds to the same function $ g_k(x)$ is attained in the case of opposite parity.
  
  \vskip+0.3cm
  \textbf{Remark 2.} Proposition 1 together with Remark 1 
calculate right hand side function $g_k(x)$ in 
  (\ref{sec12}) just for values $ \alpha =\alpha_k^{\pm 1}$. It happens that these functions are not the optimal ones.
  The optimal functions $f_k(x)$ correspond to the numbers defined in (\ref{betamain}) which are equivalent to $\alpha_k$ but have non-trivial pre-period.

\section{Proof of Theorem 2} \label{5.}
  
Now, let us  generalize the result of Proposition 1 to all the numbers equivalent to $\alpha_k$.

\subsection{The beginning}

Let $\alpha = [a_{0}; a_{1}, \,\cdots, a_{s}, \overline {1_{2k-1}, 2}]$ while $s>0$, and $a_0, a_1, \,\cdots, a_s$ be the pre-period of minimal length. 
Now as it was explained in \cite{Lesca} we can use a formula 
$$
\limsup\limits_{n \to \infty} q_{n+1}||q_{n} \alpha_{k}|| = \lim\limits_{m \to \infty}  q_{2km+s-1}||q_{2km+s-2} \alpha_{k}|| =\lim\limits_{m \to \infty}  q_{2km+s-2}||q_{2km+s-3} \alpha_{k}||
$$
similar to (\ref{precise})
and  consider subsequences $n_{m} = 2km+s-3$ or $n_{m} = 2km+s-2$.

Let $n = s + 2km - 2$ for sufficiently large $m$. In this case, we obtain
$$q_{n+1}||q_{n} \alpha|| = \frac{1}{1 + \frac{q_n}{q_{n+1}} \cdot [0; \overline {2, 1_{2k-1}}]}.$$

Similarly, we get 
$$\lim\limits_{n \to \infty} \frac{q_n}{q_{n+1}} =  \lim\limits_{n \to \infty} [0; 1_{2k-1}, \underbrace{2, 1_{2k-1}, \,\cdots, 2, 1_{2k-1}}_{m-1 \text{times}}, a_s, \,\cdots, a_1] = [0; \overline{ 1_{2k-1}, 2}] = \frac{1}{\alpha_{k}}$$.

Now we consider the finite fraction  $r_s = [a_s; a_{s-1}, \,\cdots, a_1]$ and $\gamma $
defined in (\ref{re2}).
Again we write the difference $\frac{1}{\alpha_k} - \frac{q_n}{q_{n+1}}$ as 
$$\frac{1}{\alpha_k} - \frac{q_n}{q_{n+1}} = \frac{p_{n+1-s}^{*}  \gamma + p_{n-s}^{*}}{q_{n+1-s}^{*} \gamma  + q_{n-s}^{*}} - \frac{p_{n+1-s}^{*} r_s + p_{n-s}^{*}}{q_{n+1-s}^{*} r_s + q_{n-s}^{*}}$$
(here $p_{i}^{*}$ and $q_{i}^{*}$ are, respectively, the numerator and the denominator of the $i$-th convergent of the number $\frac{1}{\alpha_k}$) and so
\begin{equation}\label{wer10}
\frac{1}{\alpha_k} - \frac{q_n}{q_{n+1}} = \frac{\gamma - r_s}{(q_{n+1-s}^{*} \gamma + q_{n-s}^{*})(q_{n+1-s}^{*} r_s + q_{n-s}^{*})}.
\end{equation}

Therefore,
\begin{equation}\label{wer1}
q_{n+1}||q_{n} \alpha|| = \frac{1}{1 + \frac{1}{\gamma}( \frac{1}{\alpha_k} - \frac{\gamma - r_s}{(q_{n+1-s}^{*} \gamma + q_{n-s}^{*})(q_{n+1-s}^{*} r_s + q_{n-s}^{*})})}.
\end{equation}

Now, let $n = s + 2km - 1$ for sufficiently large $m$. Analogously,
\begin{equation}\label{wer2}
q_{n+1}||q_{n} \alpha|| = \frac{1}{1 + \frac{1}{\alpha_k}( \frac{1}{\gamma} - \frac{r_s - \gamma}{(q_{n+1-s}^{'} \gamma + q_{n-s}^{'})(q_{n+1-s}^{'} r_s + q_{n-s}^{'})})}
,
\end{equation}
where $p_{i}^{'}$ and $q_{i}^{'}$ are  the numerator and the denominator of the $i$-th convergent of the number $\frac{1}{\gamma}$ respectively.

Note
that in one of the formulas (\ref{wer1}) and (\ref{wer2}) we have 
\begin{equation}\label{wer3}
q_{n+1}||q_{n} \alpha|| > D_k.
\end{equation}
This depends on the sign of the difference $r_s - \gamma$. Our aim is to find the minimal value of
the right hand sides of (\ref{wer1}) and (\ref{wer2}) under the condition that the sign is chosen to satisfy
(\ref{wer3}).
  This  value is minimal as a function of $q_{n+1}$ for fixed $k$ and $m$ when $|r_s - \gamma|$ is minimal.

\subsection{Inequalities with continued fractions}

We formulate here several easy auxiliary statements.

\textbf{Lemma 1.} \textit{Let $\alpha = [a_{0}; a_{1}, \,\cdots, a_{m}, \,\cdots]$, $b = [b_{0};b_{1}, \,\cdots, b_{n}]$, where $n$ is a fixed even number and $b< \alpha.$ Then the difference $\alpha - b$ is minimal when $b = [a_{0};a_{1}, \cdots, a_{n}].$}

\textbf{Proof.} For $n=0$ the statement is obvious. Assume that this is true for ${n = k}$ and prove it for ${n = k + 2}$.
So we may suppose that 
for $[a_{2}; a_{3}, \,\cdots, a_{m}, \,\cdots]$
the minimal value of
$$
[a_{2}; a_{3}, \,\cdots, a_{m}, \,\cdots]-[b_{2};b_{3}, \,\cdots, b_{n}]
$$
attains for $b = [a_{2};a_{3}, \cdots, a_{n}].$

Let in this case the minimum of the difference be achieved at the number $b = [b_0; b_1, r_k],$ where $r_k = [b_2 ; b_{3}, \,\cdots, b_{k+2}]$. Obviously, ${b_0 = a_0}$. If ${b_1 < a_1}$, then the difference ${\alpha - b}$ is negative,
so the assertion of Lemma 1 is not satisfied.
If ${b_1 > a_1}$, then $b < [a_{0}; a_{1}, \,\cdots, a_{k+2}]$, and this means  $\alpha - b > \alpha - [a_{0}; a_{1}, \cdots, a_{k+2}]>0$. 
This cannot be true as we have assumed that $\alpha - b$ is minimal.
  Therefore, $b = [a_{0}; a_{1}, r_{k}]$. However, then 
\begin{equation}\label{wer4}
\alpha - b = \frac{[a_2; a_{3}, \,\cdots] - r_{k}}{(a_{1} \cdot [a_2; a_{3}, \,\cdots] + 1)(a_{1} \cdot r_k + 1)}.\end{equation}
Right hand side of (\ref{wer4}) is a decreasing function in $r_k$.
So
  by the induction hypothesis described in the beginning of the proof, right hand side of (\ref{wer4}) is minimal when $r_{k} = [a_{2}; a_{3}, \,\cdots, a_{k+2}].$ This is just what we need.$\Box$

\vskip+0.3cm
Lemmas 2 -- 4  are similar to Lemma 1. They correspond to the cases of different parities. They can be  proven in the absolutely similar way. We left the proofs to the reader.
\vskip+0.3cm
\textbf{Lemma 2.} \textit{Let $\alpha = [a_{0}; a_{1}, \,\cdots, a_{m}, \,\cdots]$, $b = [b_{0};b_{1}, \,\cdots, b_{n}]$, where $n$ is a fixed odd number and $b< \alpha$. Then the difference $\alpha - b$ is minimal when $b = [a_{0};a_{1}, \,\cdots, a_{n-1}, a_{n} + 1].$}

\vskip+0.3cm
\textbf{Lemma 3.} \textit{Let $\alpha = [a_{0}; a_{1}, \,\cdots, a_{m}, \,\cdots]$, $b = [b_{0};b_{1}, \,\cdots, b_{n}]$, where $n$ is a fixed odd number and $b > \alpha.$ Then
n the difference $b - \alpha$ is minimal when $b = [a_{0};a_{1}, \,\cdots, a_{n-1}, a_{n}].$}
\vskip+0.3cm
\textbf{Lemma 4.} \textit{Let $\alpha = [a_{0}; a_{1}, \,\cdots, a_{m}, \,\cdots]$, $b = [b_{0};b_{1}, \,\cdots, b_{n}]$, where $n$ is a fixed even number and $b > \alpha.$ Then the difference $b - \alpha$ is minimal when $b = [a_{0};a_{1}, \,\cdots, a_{n-1}, a_{n} + 1].$}

\vskip+0.3cm
The next lemma is a little bit more complicated.
\vskip+0.3cm
\textbf{Lemma 5.} \textit{ 
Let us fix $k\ge 1$ and 
consider the sum of continued fractions
$$
W_m(s)=[0;1_{s-1}, 2,\underbrace{1_{2k-1}, 2, \,\cdots\,, 1_{2k-1}, 2}_m] + [ 1; 1_{2k-s-1}, \overline {2, 1_{2k-1}}]$$ 
for different values of $s$ in the interval $1\le s \le 2k$.
 Then for $m$ large enough,
$$
\max_{1\le s\le 2k,\,\, s\,\, \text{is odd}}
W_m(s) = 
\begin{cases}
    W_m(k) \,\,\, \text{if}\,\, $k$ \,\,\text{is odd},\cr
    W_m(k-1) \,\,\text{or}\,\, W_m(k+1) \,\,\, \text{if}\,\, $k$ \,\,\text{is even};
\end{cases}
$$
and
$$
\max_{1\le s\le 2k,\,\, s\,\, \text{is even}}
W_m(s) = 
[0;\underbrace{1_{2k-1}, 2, \,\cdots, \,1_{2k-1}, 2}_{m+1}] + [ \overline {2, 1_{2k-1}}],
$$
and the right hand side in the last formula may be considered as $W_m(2k)$.
}

\vskip+0.3cm

\textbf{Proof}.
First of all we consider the value
$$
W_\infty (s) = [0;1_{s-1}, \overline{ 2, 1_{2k-1}}] + [0; 1_{2k-s-1}, \overline {2, 1_{2k-1}}].
$$
We consider differences
$$[0;1_{s-1}, \overline{2, 1_{2k-1}}] - [0;1_{s-3}, \overline{2,1_{2k-1}}] = \frac{(-1)^{s-3}([\overline{2, 1_{2k-1}}] - [1;1, \overline{2, 1_{2k-1}}])}{(F_{s-2} [\overline{2, 1_{2k-1}}] + F_{s-3})(F_{s-2}[1;1, \overline{2, 1_{2k-1}}] + F_{s-3})}$$
 and
 \begin{multline*}
[0;1_{2k-s+1}, \overline{2, 1_{2k-1}}] - [0;1_{2k-s-1}, \overline{2,1_{2k-1}}] =\\
=\frac{(-1)^{2k-s-1}([\overline{2, 1_{2k-1}}] - [1;1, \overline{2, 1_{2k-1}}])}{(F_{2k-s} [\overline{2, 1_{2k-1}}] + F_{2k-s-1})(F_{2k-s}[1;1, \overline{2, 1_{2k-1}}] + F_{2k-s-1})},
\end{multline*}
where 
$F_i$ stands for $i$-th  Fibonacci number.
\vskip+0.3cm

Let $s$ be odd and then $ (-1)^{s-3}=(-1)^{2k-s-1}=1$.

 We see that for $s-1 > 2k-s-1$ 
 one has
$$[0;1_{s-1}, \overline{2, 1_{2k-1}}] - [0;1_{s-3}, \overline{2,1_{2k-1}}] \leq [0;1_{2k-s+1}, \overline{2, 1_{2k-1}}] - [0;1_{2k-s-1}, \overline{2,1_{2k-1}}],$$
and this is  equivalent to
$$
W_\infty(s)=[0;1_{s-1}, \overline{2, 1_{2k-1}}] + [0;1_{2k-s-1}, \overline{2,1_{2k-1}}] \leq 
$$
\begin{equation}\label{bub1}
\leq [0;1_{s-3}, \overline{2,1_{2k-1}}] + [0;1_{2k-s+1}, \overline{2, 1_{2k-1}}]=W_\infty (s-2).
\end{equation}
The equality in (\ref{bub1}) happens
 only when $s-3 = 2k-s-1$, that is, when $s = k+1$ (the fact that $s$ is odd means that in this case $k$ has to be even). 

Similarly when ${s-1 < 2k-s-1}$, we obtain that 
$$
W_{\infty}(s)= [0;1_{s-1}, \overline{2, 1_{2k-1}}] + [0;1_{2k-s-1}, \overline{2,1_{2k-1}}] \leq $$
\begin{equation}\label{bub2}[0;1_{s+1}, \overline{2,1_{2k-1}}] + [0;1_{2k-s-3}, \overline{2, 1_{2k-1}}] = W_{\infty}(s+2).
\end{equation}
The equality in (\ref{bub2}) happens only when $s-1 = 2k-s-3$, that is, when $s = k-1$ (again the fact that $s$ is odd means that $k$ has to be even).
So in the case when $k$ is even, we get 
$$
\max_{1\le s\le 2k,\,\, s\,\, \text{is odd}}
W_\infty(s) =  W_{\infty}(k-1) = W_{\infty}(k+1),
$$
meanwhile for odd $k$ the maximum is attained just at the unique value of $s$:
$$
\max_{1\le s\le 2k,\,\, s\,\, \text{is odd}} =
 W_{\infty}(k).  
 $$
So we conclude that 
$$
\max_{1\le s\le 2k,\,\, s\,\, \text{is odd}}
W_\infty(s)= \begin{cases}
    W_{\infty}(k)  \,\,\, \text{if}\,\, $k$ \,\,\text{is odd},\cr
    W_{\infty}(k-1) = W_{\infty}(k+1) \,\,\, \text{if}\,\, $k$ \,\,\text{is even};
\end{cases}
$$

 \vskip+0.3cm
The case of even $s$ is treated in a completely similar way, with the difference that inequalities (\ref{bub1}), (\ref{bub2}) are replaced by their opposites. Then the maximum will be attained at the endpoint of the interval of admissible values of $s$ and so
$$\max_{1\le s\le 2k,\,\, s\,\, \text{is even}}
W_\infty(s) = W_{\infty}(2k).$$

 \vskip+0.3cm
 
Now we deal with the quantity  $W_{m}(s)$ 
from the formulation of  Lemma 5 for a finite value of $m$.
We show that for all $m$ large enough it attains its maximum at the same $s$ as
 $W_{\infty}(s)$. 
Indeed, 
let for infinite set of values of
  $m$ one has $$\max_{1\le s\le 2k}
W_m(s) = W_{m}(s_0),$$
but
 $$\max_{1\le s\le 2k}
W_\infty(s) = W_{\infty}(s_1) \neq W_{\infty}(s_0).$$
We should note that for any fixed $s$ one has 
$$\lim\limits_{m \to \infty} W_{m}(s) = W_{\infty}(s)+1.
$$
But then
there exists 
  $m'$
  such that  for all 
  $m>m'$ we have $W_{m}(s_1) \neq W_{m}(s_0)$. This is a contradiction. 
Lemma 5 is proven.$\Box$




\subsection{Final analysis of cases}

Here we consider four cases. In each of the cases we find the approximations satisfying
(\ref{wer3}). Among them we choose the minimal one.

\vskip+0.3cm
\textbf{1. $n = 2mk + s - 2$, $s$ is odd.}

In this case, the value of $q_{n+1}||q_{n} \alpha||$ is minimal while still greater than $D_k$ when
the difference
$$\gamma - r_s = [\overline{2, 1_{2k-1}}] - [a_{s}; a_{s-1}, \,\cdots, a_{1}] = [b_0; b_1, \,\cdots] - [a_{s}; a_{s-1}, \,\cdots, a_1],\,\,\,
\gamma = [b_0;b_1,...]
$$
attains its minimal positive value. Since the number $s-1$ is even, by Lemma 1 we obtain that the difference under consideration takes its minimal positive value in the case where $a_s = b_0, \  a_{s-1} = b_1, \  \cdots, \  a_{2} = b_{s-2}, \  a_1 = b_{s-1} $. So, if $s \geqslant 2k+1$,
we can write
$$
s = 2kv+l
\,\,\,\text{and}\,\,\,
r_s = [\underbrace{2; 1_{2k-1} \,\cdots, 2,1_{2k-1}}_{v \  \text{times}}, 2, 1_{l-1}].
$$ 
But
we have assumed that our 
  pre-period 
  has the minimal length. So $v=0$.
  Thus, we can assume that $s=l \leqslant 2k-1$ and $r_s = [2; 1_{s-1}],$
  or in other words 
  $$ \alpha = [0;1_{s-1},\overline {2, 1_{2k-1}}]
  =[0;1_{s-1},2, \alpha_k]
  .$$

Then simplification of formula (\ref{wer1}) gives us
\begin{equation}\label{kk1}
q_{n+1} ||q_{n} \alpha|| = 
q_{2mk + s - 1} ||q_{2mk + s - 2} \alpha|| 
=
\frac{1}{1+ \frac{1}{\gamma}(\frac{1}{\alpha_{k}} - \frac{1}{q_{n+1}^{2}([0;1_{s-1}, 2, \,\cdots, 1_{2k-1}, 2] + [ 1; 1_{2k-s-1}, \overline {2, 1_{2k-1}}])})}.
\end{equation}
We give the following explanation.
Equality (\ref{wer1}) is obtained from (\ref{wer10}) where 
$\frac{q_n}{q_{n+1}}$ is not necessarily a convergent fraction to 
$\frac{1}{\alpha_k}$. But in the case under consideration
$$\frac{q_n}{q_{n+1}}
=[0; \underbrace{1_{2k-1},2, \cdots, 1_{2k-1}, 2}_{m \  \text{times}}, 1_{s-1}]
$$ is a convergent fraction to $\frac{1}{\alpha_k}$,
and everything is simpler.

 The quantity in (\ref{kk1}) is minimal as a function of $q_{n+1}$ for fixed $k$ and $m$ when the value of $[0;1_{s-1}, 2, \,\cdots, 1_{2k-1}, 2] + [ 1; 1_{2k-s-1}, \overline {2, 1_{2k-1}}]$ is maximal. By Lemma 5, this is achieved when $s = k$ for odd $k$ and when $s = k+1$ or $s = k-1$ for even $k$. 

\vskip+0.3cm
\textbf{2. $n = 2mk + s - 2$, $s$ is even.}

By Lemma 2 in this case we have $r_s = [2;1_{s-2}, 2]$, where $s \leqslant 2k.$ However, since $r_{s} = [2;1_{s-2}, 2] = [2;1_{s}]$, the situation is reduced to the previous case, as we increase $s$ by 1 and this changes the parity.

\vskip+0.3cm
\textbf{3. $n = 2mk + s - 1$, $s$ is even.}

Similarly, by Lemma 3, $r_s = [2;1_{s-1}]$ for $s \leqslant 2k$. Then let us rewrite $q_{n+1} ||q_{n} \alpha||$ as
$$q_{n+1} ||q_{n} \alpha|| = \frac{1}{1+ \frac{1}{\alpha_k}(\frac{1}{\gamma} - \frac{1}{q_{n+1}^{2}([0;1_{s-1}, 2, \,\cdots, 1_{2k-1}, 2] + [ 1; 1_{2k-s-1}, \overline {2, 1_{2k-1}}])})}.$$

This expression is minimal as a function of $q_{n+1}$ for fixed $k$ and $m$ when the value of $[0;1_{s-1}, 2, \,\cdots, 1_{2k-1}, 2] + [ 1; 1_{2k-s-1}, \overline {2, 1_{2k-1}}]$ is maximal. By Lemma 5, this is achieved when $s=2k$, but this is equivalent to the case where there is no pre-period at all.

\vskip+0.3cm
\textbf{4. $n = 2mk + s - 1$, $s$ is odd.}

Similarly, by Lemma 4, $r_s = [2;1_{s-2},2]$ for $s \leqslant 2k-1$, but since $r_{s} =  [2;1_{s-2},2] = [2;1_{s}]$, this situation reduces to the previous case {\bf 3}.

\vskip+0.3cm
Consideration of cases {\bf 1} - {\bf 4} show that 
  the  numbers 

\vskip+0.3cm
  \noindent
  $
  \beta_{k} = [0;1_{k-1}, \overline{2, 1_{2k-1}}]$ when $k$ is odd and 
  
    \noindent
  $\beta_{k}^{(1)} = [0;1_{k}, \overline{2, 1_{2k-1}}]$ or $\beta_{k}^{(2)} = [0;1_{k-2}, \overline{2, 1_{2k-1}}]$ when $k$ is even

\vskip+0.3cm
\noindent
from Theorem 2
have optimal pre-periods, and we simply must calculate optimal approximation functions by means of these numbers.

Calculation show that 
  for the number $\beta_{k}$ when $k$ is odd. 
  we have
\begin{equation}\label{fin1}
c = \frac{2D_{k}}{1 + \sqrt{ 1 - \frac{2 \alpha_{k}}{(2 \beta_{k} + 1)( \alpha_{k} + 1) q_{n+1}^{2}}}} = f_k(q_{n+1})\,\,\,
\text{for}\,\,\,
 n = (2m+1)k - 2.
 \end{equation}

Similarly, for the number $\beta_{k}^{(1)}$ or $\beta_{k}^{(2)}$ when $k$ is even, we have 
\begin{equation}\label{fin2}
q_{n+1}||q_{n} \beta_{k}^{(i)}|| = \frac{2D_{k}}{1 + \sqrt{ 1 - \frac{2 \alpha_{k}}{(\beta^{(1)}_{k} + \beta^{(2)}_{k} + 1)( \alpha_{k} + 1) q_{n+1}^{2}}}} = f_k(q_{n+1})
\end{equation}
$$
\text{for}\,\,\,
\begin{cases}
    n = (2m+1)k - 1  \,\,\,\text{for}\,\,\,i = 1,\cr
    n = (2m+1)k - 3\,\,\,\text{for}\,\,\,i = 2.
\end{cases}
$$

\vskip+0.3cm

Finally, Statement 1 of Theorem 2 follows from our construction, because 
numbers $\beta= \beta_{k}, \beta_{k}^{(1)}, \beta_{k}^{(2)}\sim \alpha_k$ are chosen 
in such a way, that the values $q_{n+1}||q_n\beta||$ with the optimal choice of $n$  satisfy  $q_{n+1}||q_n\beta||>D_k$ and are minimal among all numbers equivalent to $\alpha_k$.

As for Statement 2,  equalities (\ref{fin1}) and (\ref{fin2}) ensure property b), while property a) follows from the fact that for the values of n
which do not satisfy conditions from 
(\ref{fin1},\ref{fin2}) we have 
$q_{n+1}||q_{n} \beta ||<D_k$.

 Theorem 2 is proven.$\Box$

\end{document}